\documentclass{amsart}
\usepackage[utf8]{inputenc}
\usepackage[T1]{fontenc}
\usepackage{lmodern}
\usepackage{amsmath,amssymb,amsthm,mathtools}
\usepackage{booktabs,array}
\usepackage{tikz-cd}
\allowdisplaybreaks
\usetikzlibrary{arrows.meta}
\tikzcdset{bighook/.style={{Hooks[right,length=1.1ex,width=2.6ex,line width=0.7pt]}->}}
\usepackage{microtype}
\usepackage[hidelinks]{hyperref}

\newtheorem{theorem}{Theorem}[section]
\newtheorem{proposition}[theorem]{Proposition}

\theoremstyle{definition}
\newtheorem{definition}[theorem]{Definition}
\newtheorem{example}[theorem]{Example}
\newtheorem{problem}[theorem]{Problem}
\theoremstyle{remark}
\newtheorem{remark}[theorem]{Remark}

\newcommand{\C}{\mathbb{C}}
\newcommand{\K}{\mathbb{K}}
\newcommand{\lt}{\dashv}
\newcommand{\rt}{\vdash}
\newcommand{\Ann}{\operatorname{Ann}}
\newcommand{\Aut}{\operatorname{Aut}}
\newcommand{\GL}{\mathrm{GL}}
\newcommand{\cat}[1]{\mathbf{#1}}

\begin{document}

\title[Classification of diassociative and dendriform algebras]{Classification problems for diassociative and dendriform algebras: a survey}

\author{Ikrom M. Rikhsiboev}

\subjclass[2020]{Primary 17-02; Secondary 17A30, 17A32, 16W99}
\keywords{diassociative algebra, dendriform algebra, Zinbiel algebra, classification, central extension, degeneration}

\maketitle

\begin{center}\small
Westminster International University in Tashkent, Tashkent, Uzbekistan\\
E-mail: \href{mailto:irikhsiboev@wiut.uz}{irikhsiboev@wiut.uz}
\end{center}

\begin{abstract}
Diassociative algebras (associative dialgebras) and dendriform algebras are vector spaces with two binary operations, introduced by J.-L.~Loday in the study of Leibniz algebras; their operads are Koszul dual to each other. We survey the classification problem for these two classes, including the Zinbiel algebras, the commutative dendriform algebras, which are Koszul dual to Leibniz algebras. We recall the definitions, the main examples and the Loday diagram with its extension by left-symmetric algebras and dialgebras, and describe the methods used to classify algebras of small dimension: direct computation from lists of associative algebras, annihilator (central) extensions and geometric classification by degenerations. We then collect, by dimension, the known algebraic classifications of diassociative algebras, dendriform algebras and their commutative case, the Zinbiel algebras, point out lists that later work corrected, and summarise the cohomological results that the methods rely on. We close with a list of open problems.
\end{abstract}

\section{Introduction}

The structure theory of associative algebras goes back to the second half of the nineteenth century. Peirce's \emph{Linear associative algebra} \cite{Peirce1881} already contained lists of algebras of small dimension, and the classification of associative algebras of a given dimension remains a source of problems today \cite{Pierce1982,RRB09,deGraaf18}. In the last decades of the twentieth century non-associative structures moved to the centre of algebra, largely because of their role in geometry and physics. Lie algebras are the main example: over a field of characteristic zero every finite-dimensional Lie algebra is the semidirect sum of a semisimple subalgebra and its solvable radical, the semisimple part is a direct sum of simple algebras classified by Cartan and Killing, and the problem reduces to solvable and finally to nilpotent algebras \cite{Jacobson,GozeKhakimdjanov,BFNT}.

In the early 1990s J.-L.~Loday introduced a family of algebraic structures, now often called \emph{Loday algebras}, while studying the periodicity of algebraic $K$-theory and Leibniz homology \cite{Loday93,LodayEtAl}. The four main members are Leibniz algebras, diassociative algebras, dendriform algebras and Zinbiel algebras. A Leibniz algebra is a non-antisymmetric analogue of a Lie algebra. A \emph{diassociative algebra} (associative dialgebra) has two associative products $\lt$ and $\rt$ and plays for Leibniz algebras the role that associative algebras play for Lie algebras: the bracket $[x,y]=x\lt y-y\rt x$ is a Leibniz bracket. A \emph{dendriform algebra} has two products $\prec$ and $\succ$ whose sum is associative, so it ``splits'' an associative product; Zinbiel algebras are the commutative dendriform algebras. In the language of operads, $\cat{Dias}$ and $\cat{Dend}$ are Koszul dual, as are $\cat{Leib}$ and $\cat{Zinb}$ \cite{GinzburgKapranov,LodayEtAl}. These structures have since appeared in homological algebra \cite{Frabetti97,Frabetti98,Frabetti01}, deformation theory \cite{MajumdarMukherjee03,MajumdarMukherjee05,Yau08,FialowskiMajumdar}, combinatorics of planar trees and polytopes \cite{LodayRonco04}, Rota--Baxter algebras \cite{Aguiar00,EbrahimiFard02,EbrahimiFardGuo08}, polynomial identities \cite{Kolesnikov08,Bremner12,KolesnikovVoronin} and low-dimensional topology \cite{KravchenkoPolyak}.

This paper surveys one problem for the two classes with two binary operations, diassociative and dendriform algebras: their \emph{classification}. As for associative and Lie algebras, a complete classification in all dimensions is out of reach, and the literature consists of lists of algebras of small dimension, often restricted to nilpotent algebras, together with methods to obtain them. These results are scattered over journals and conference proceedings, several lists were later corrected, and the methods developed for other varieties (central extensions, degenerations) have been applied to dialgebras only partly. Our aim is to collect the known classifications by dimension, explain the methods behind them, record the corrections and list the problems that remain open. Earlier surveys covering related ground are Bremner's survey of dialgebras and polynomial identities \cite{Bremner12}, Rakhimov's survey of the classification problem for Loday algebras \cite{Rakhimov16b}, and the surveys of Kaygorodov and coauthors on the algebraic and geometric classification of non-associative algebras \cite{Kaygorodov24,KKP24}.

Section~\ref{sec:def} recalls the definitions, the main examples and the Loday diagram with its extension by left-symmetric algebras and dialgebras. Section~\ref{sec:methods} describes the classification problem and its methods. Sections~\ref{sec:dias} and~\ref{sec:dend} collect the results for diassociative and dendriform algebras, Section~\ref{sec:coh} summarises the cohomological results used by the methods, and Section~\ref{sec:open} lists open problems.

\section{Algebras with two binary operations}\label{sec:def}

Throughout, an \emph{algebra} is a vector space over a field $\K$ with one or more bilinear products; the dimension of the algebra is the dimension of the space. Unless stated otherwise, $\K=\C$.

\subsection{Leibniz algebras}
A (right) \emph{Leibniz algebra} is a vector space $L$ with a bilinear bracket satisfying
\begin{equation}\label{eq:leib}
[[x,y],z]=[[x,z],y]+[x,[y,z]]\qquad(x,y,z\in L).
\end{equation}
An antisymmetric Leibniz algebra is a Lie algebra. Leibniz algebras were introduced by Bloh \cite{Bloh} and, independently, by Loday \cite{Loday93}. Their structure theory is well developed (Levi decomposition, Cartan subalgebras, solvable and nilpotent algebras), and their low-dimensional classification is a large subject in its own right; see \cite{AyupovOmirov98,AyupovOmirov99,AyupovOmirov01,AOR06} and the survey \cite{Kaygorodov24}. In this paper Leibniz algebras appear mainly as the target of the functor from diassociative algebras.

\subsection{Diassociative algebras}

\begin{definition}[\cite{LodayEtAl}]
A \emph{diassociative algebra} (or \emph{associative dialgebra}) is a vector space $D$ with two bilinear products $\lt$ and $\rt$ such that, for all $x,y,z\in D$,
\begin{equation}\label{eq:dias}
\begin{aligned}
(x\lt y)\lt z&=x\lt(y\lt z), & (x\lt y)\lt z&=x\lt(y\rt z), & (x\rt y)\lt z&=x\rt(y\lt z),\\
(x\lt y)\rt z&=x\rt(y\rt z), & (x\rt y)\rt z&=x\rt(y\rt z).&&
\end{aligned}
\end{equation}
\end{definition}

The first and the last identity say that $\lt$ and $\rt$ are associative, and the three middle identities link the two products. Many authors therefore define a diassociative algebra as a vector space with two associative products satisfying the three identities
\[
\begin{gathered}
(x\lt y)\lt z=x\lt(y\rt z),\qquad (x\rt y)\lt z=x\rt(y\lt z),\\
(x\lt y)\rt z=(x\rt y)\rt z,
\end{gathered}
\]
the last one being the combination of the fourth and fifth identities in \eqref{eq:dias}; the two definitions are equivalent. An associative algebra is a diassociative algebra with $x\lt y=x\rt y$. The following examples are due to Loday \cite{LodayEtAl} and Lin--Zhang \cite{LinZhang}.

\begin{example}
Let $(A,d)$ be a differential associative algebra with $d^2=0$. Then $x\lt y=x\,d(y)$, $x\rt y=d(x)\,y$ define a diassociative algebra on $A$.
\end{example}

\begin{example}
Let $A$ be an associative algebra, $M$ an $A$-bimodule and $f\colon M\to A$ a bimodule map. Then $m\lt m'=m\,f(m')$, $m\rt m'=f(m)\,m'$ define a diassociative algebra on $M$.
\end{example}

\begin{example}
Let $A$ be an associative algebra and $D=A^n$. Then $(x\lt y)_i=x_i\bigl(\sum_j y_j\bigr)$ and $(x\rt y)_i=\bigl(\sum_j x_j\bigr)y_i$ define a diassociative algebra.
\end{example}

\begin{example}[\cite{LinZhang}]
On the polynomial algebra $\K[x,y]$ the products \[f(x,y)\lt g(x,y)=f(x,y)g(y,y),\qquad f(x,y)\rt g(x,y)=f(x,x)g(x,y)\] define a diassociative algebra.
\end{example}

\begin{proposition}[\cite{LodayEtAl}]
If $D$ is diassociative, then $[x,y]=x\lt y-y\rt x$ satisfies \eqref{eq:leib}. This gives a functor $\cat{Dias}\to\cat{Leib}$ extending the commutator functor $\cat{As}\to\cat{Lie}$.
\end{proposition}

Loday showed that the universal enveloping algebra of a Leibniz algebra is naturally a diassociative algebra, and Aymon and Grivel proved a Poincar\'e--Birkhoff--Witt theorem in this setting \cite{AymonGrivel}. The free diassociative algebra has a basis described by planar binary trees \cite{LodayEtAl}.

\subsection{Dendriform and Zinbiel algebras}

\begin{definition}[\cite{LodayEtAl}]
A \emph{dendriform algebra} is a vector space $E$ with two bilinear products $\prec$ and $\succ$ such that, for all $x,y,z\in E$,
\begin{equation}\label{eq:dend}
\begin{aligned}
(x\prec y)\prec z&=x\prec(y\prec z+y\succ z),\qquad (x\succ y)\prec z=x\succ(y\prec z),\\
(x\prec y+x\succ y)\succ z&=x\succ(y\succ z).
\end{aligned}
\end{equation}
\end{definition}

Adding the three identities shows that $x*y=x\prec y+x\succ y$ is associative, which gives a functor $\cat{Dend}\to\cat{As}$. Moreover $x\circ y=x\succ y-y\prec x$ is a left-symmetric (pre-Lie) product, which gives a functor $\cat{Dend}\to\cat{LSA}$ \cite{Aguiar00,LodayEtAl}.

\begin{example}[\cite{Aguiar00,EbrahimiFard02}]
Let $A$ be an associative algebra with a Rota--Baxter operator $R$ of weight $0$, that is, $R(x)R(y)=R\bigl(R(x)y+xR(y)\bigr)$. Then $x\prec y=xR(y)$ and $x\succ y=R(x)y$ define a dendriform algebra. Operators of non-zero weight give dendriform trialgebras in the same way \cite{EbrahimiFard02}, and every dendriform algebra embeds into a Rota--Baxter algebra \cite{EbrahimiFardGuo08,ChenMo,GubarevKolesnikov}.
\end{example}

\begin{definition}
A \emph{Zinbiel algebra} is a vector space $Z$ with a product satisfying $(x\cdot y)\cdot z=x\cdot(y\cdot z)+x\cdot(z\cdot y)$.
\end{definition}

A Zinbiel algebra is exactly a dendriform algebra with $x\succ y=y\prec x$, through $x\cdot y=x\prec y$; its symmetrisation $x\cdot y+y\cdot x$ is commutative and associative. A standard example is the space of integrable functions with $(f\cdot g)(x)=f(x)\int_0^xg(t)\,dt$. The name is Leibniz read backwards, and Loday also wrote under the pseudonym G.~W.~Zinbiel \cite{Zinbiel}.

\subsection{The Loday diagram}
\begin{figure}[ht]
\[
\begin{tikzcd}[column sep=2.1cm,row sep=1.3cm]
\cat{Zinb} \arrow[r, bighook] \arrow[d, "{x\cdot y+y\cdot x}"'] & \cat{Dend} \arrow[d, "{\prec+\succ}"] & \\
\cat{Com} \arrow[r, bighook] & \cat{As} \arrow[r, bighook] \arrow[d, "{xy-yx}"'] & \cat{Dias} \arrow[d, "{x\lt y-y\rt x}"]\\
& \cat{Lie} \arrow[r, bighook] & \cat{Leib}
\end{tikzcd}
\]
\caption{The Loday diagram. Hooked arrows are inclusions of subcategories; the other arrows are the functors described in Section~\ref{sec:def}.}\label{fig:loday}
\end{figure}
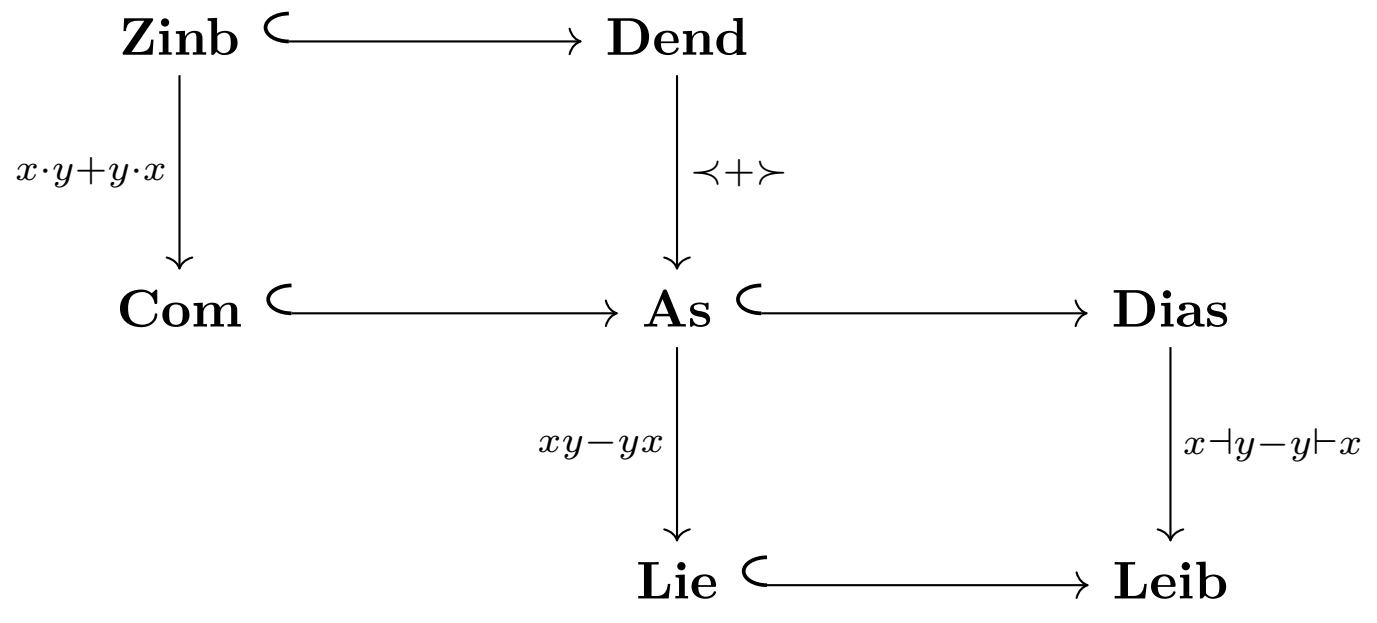

The operads $\cat{Dias}$ and $\cat{Dend}$ are binary, quadratic and Koszul, and $\cat{Dias}^!=\cat{Dend}$; likewise $\cat{Leib}^!=\cat{Zinb}$, while $\cat{As}^!=\cat{As}$ and $\cat{Com}^!=\cat{Lie}$ \cite{GinzburgKapranov,LodayEtAl}. Chapoton \cite{Chapoton05} showed that these operads are anticyclic and computed their characters. The corresponding categories of algebras form the commutative diagram of functors in Figure~\ref{fig:loday}, in which the Koszul dual categories occupy symmetric positions \cite{LodayEtAl}.

\subsection{The extended Loday diagram}
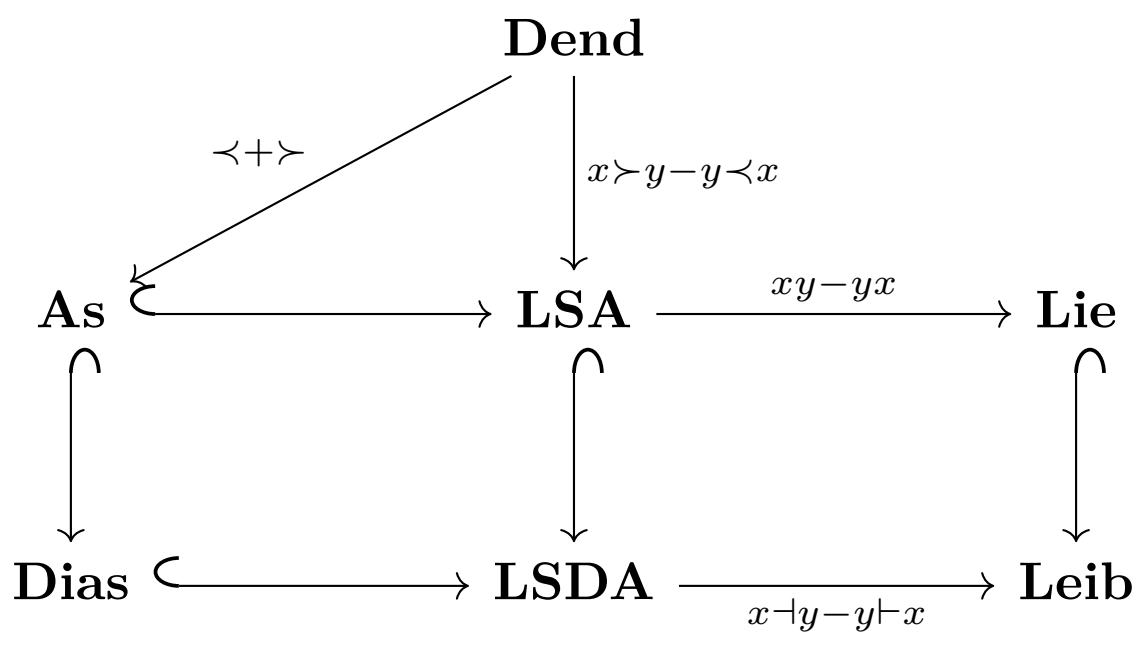
\begin{figure}[ht]
\[
\begin{tikzcd}[column sep=2.1cm,row sep=1.3cm]
& \cat{Dend} \arrow[dl, "{\prec+\succ}"'] \arrow[d, "{x\succ y-y\prec x}"] & \\
\cat{As} \arrow[r, bighook] \arrow[d, bighook] & \cat{LSA} \arrow[r, "{xy-yx}"] \arrow[d, bighook] & \cat{Lie} \arrow[d, bighook]\\
\cat{Dias} \arrow[r, bighook] & \cat{LSDA} \arrow[r, "{x\lt y-y\rt x}"'] & \cat{Leib}
\end{tikzcd}
\]
\caption{The extended Loday diagram with left-symmetric algebras and dialgebras.}\label{fig:ext}
\end{figure}

Left-symmetric algebras (LSA) satisfy $(x,y,z)=(y,x,z)$ for the associator $(x,y,z)=(xy)z-x(yz)$; see the survey \cite{Burde}. Felipe \cite{Felipe11} introduced their dialgebra analogue, the \emph{left-symmetric dialgebras} (LSDA): vector spaces with products $\lt,\rt$ satisfying
\[
\begin{aligned}
&x\lt(y\lt z)=x\lt(y\rt z),\qquad (x\rt y)\rt z=(x\lt y)\rt z,\\
&x\lt(y\lt z)-(x\lt y)\lt z=y\rt(x\lt z)-(y\rt x)\lt z,\\
&x\rt(y\rt z)-(x\rt y)\rt z=y\rt(x\rt z)-(y\rt x)\rt z.
\end{aligned}
\]
Every diassociative algebra is an LSDA, an LSDA with $\lt=\rt$ is an LSA, and $x\lt y-y\rt x$ is again a Leibniz bracket \cite{Felipe11,Rikhsiboev16}. Adding these classes gives the extended diagram of Figure~\ref{fig:ext} \cite{Rikhsiboev16}. The two ways from $\cat{Dend}$ to $\cat{Lie}$ agree, since $(x\succ y-y\prec x)-(y\succ x-x\prec y)=x*y-y*x$.

\section{The classification problem and its methods}\label{sec:methods}

\subsection{Structure constants and orbits}
Fix a basis $e_1,\dots,e_n$. An algebra with two products is given by structure constants $e_i\lt e_j=\sum_k\gamma_{ij}^ke_k$ and $e_i\rt e_j=\sum_k\delta_{ij}^ke_k$, a point of $\C^{2n^3}$. The identities \eqref{eq:dias} (or \eqref{eq:dend}) are polynomial equations of degree two in these constants, so the $n$-dimensional diassociative algebras form an affine variety $\cat{Dias}_n$, and similarly $\cat{Dend}_n$. The group $\GL_n(\C)$ acts on the variety by change of basis, and isomorphism classes are the orbits. To \emph{classify} $n$-dimensional algebras means to give one representative of each orbit together with a criterion deciding when two representatives are isomorphic. Because the number of orbits is infinite already in small dimension, the answer usually contains parametric families. The methods and their difficulty depend on the base field and on the class of algebras considered (all, nilpotent, solvable, with given properties).

\subsection{Nilpotency}
For a diassociative algebra $D$ put $D^1=D$ and
\[D^{k}=\sum_{i+j=k}\bigl(D^i\lt D^j+D^i\rt D^j\bigr)\quad(k\ge2).\] The algebra is \emph{nilpotent} if $D^s=0$ for some $s$ \cite{Rikhsiboev07}. An analogue of Engel's theorem holds, and if the associative algebras $(D,\lt)$ and $(D,\rt)$ are nilpotent (solvable), then so is $D$ \cite{Rikhsiboev07}. The same definition applies to dendriform and other algebras with two products. Every finite-dimensional Zinbiel algebra is nilpotent: Dzhumadildaev and Tulenbaev proved this in characteristic zero \cite{DzhumadildaevTulenbaev}, and Towers extended it to all fields \cite{Towers}.

\subsection{Direct method from associative lists}
Since both products of a diassociative algebra are associative, one may start from a list of $n$-dimensional associative algebras \cite{RRB09}, fix $\lt$ to be one of them, treat the structure constants of $\rt$ as unknowns, and solve the system \eqref{eq:dias}; the solutions are then sorted into isomorphism classes. This is the approach of the first classifications of diassociative algebras in dimensions two and three \cite{BasriRikhsiboev07,RRB10,Rikhsiboev07}, carried out with computer algebra. The same idea, with the associative algebra $(E,\prec+\succ)$ fixed first, was used for dendriform algebras \cite{RRB10dend}, and in the compatible version for anti-dendriform algebras \cite{AANS}. The method is elementary, but the systems grow quickly with the dimension, and sorting the solutions into orbits without omissions is the delicate step.

\subsection{Annihilator extensions}
Every non-zero nilpotent algebra $S$ has a non-zero annihilator $\Ann(S)$, the set of elements whose products with $S$ vanish on both sides for both operations, and $S$ is an annihilator (central) extension of the smaller algebra $S/\Ann(S)$. This reduces the classification of nilpotent algebras of dimension $n$ to two steps: compute the second cohomology $H^2(A,\C)$ of each nilpotent algebra $A$ of dimension $n-s$, and find the orbits of $\Aut(A)$ on the $s$-dimensional subspaces of $H^2(A,\C)$ satisfying a non-degeneracy condition. The method is due to Skjelbred and Sund for Lie algebras \cite{SkjelbredSund}; it was developed for many varieties by Hegazi, Abdelwahab and Calder\'on Mart\'in \cite{HAC} and by Kaygorodov and his coauthors \cite{Kaygorodov24}, and for nilpotent associative algebras by de Graaf \cite{deGraaf18}. For associative dialgebras Rakhimov established the basic properties of central extensions and applied the method to low-dimensional nilpotent algebras \cite{Rakhimov16a}. Mainellis studied the related multipliers and cover algebras of nilpotent diassociative algebras \cite{Mainellis22,Mainellis23}. The method has been adapted to left-symmetric dialgebras, which contain the diassociative algebras, in \cite{RikhsiboevLSDA}.

\subsection{Characteristic sequences and gradings}
For Lie, Leibniz and Zinbiel algebras of larger dimension one classifies subclasses defined by invariants, such as the characteristic sequence (the Jordan type of generic right multiplication operators) or the natural grading. This gives the filiform and quasi-filiform classes in arbitrary dimension \cite{GozeKhakimdjanov,AdashevOmirovKhudoyberdiyev10,AdashevKhudoyberdiyevOmirov,AdashevCamachoGomezVidalKarimdjanov}. For diassociative and dendriform algebras this approach has hardly been used.

\subsection{Geometric classification}
The orbit closures in $\cat{Dias}_n$ or $\cat{Dend}_n$ define degenerations: $A$ degenerates to $B$ if $B$ lies in the closure of the orbit of $A$. The \emph{geometric classification} describes the irreducible components of the variety and the rigid algebras, whose orbits are open. It is surveyed in \cite{KKP24}. For the Loday algebras considered here it is complete for Zinbiel algebras up to dimension five \cite{KPPV18,AFK22}; for diassociative and general dendriform algebras we are not aware of any geometric classification.

\section{Diassociative algebras}\label{sec:dias}

Table~\ref{tab:dias} lists the known classifications of diassociative algebras and closely related algebras.

\begin{table}[ht]
\caption{Classification results for diassociative and related algebras.}\label{tab:dias}
\small
\begin{tabular}{@{}>{\raggedright}p{3.1cm}>{\raggedright}p{1.7cm}>{\raggedright}p{5.2cm}l@{}}
\toprule
Class & Dimension & Result & Reference\\
\midrule
Diassociative, $\C$ & 2 & complete list (four classes) & \cite{BasriRikhsiboev07}\\
Diassociative, any field & 2 & complete list; corrects and completes \cite{BasriRikhsiboev07} & \cite{Rakhimov23}\\
Diassociative, $\C$ & 3 & complete list & \cite{RRB10}\\
Nilpotent diass., $\C$ & $\le3$ & complete list; Engel-type theorem & \cite{Rikhsiboev07}\\
Nilpotent diass., $\C$ & 4 & complete list & \cite{BRR15}\\
Nilpotent diass., $\C$ & small & central extension method & \cite{Rakhimov16a}\\
Hom-associative dialgebras, $\C$ & $\le4$ & classification, $\alpha^k$-derivations, centroids & \cite{Zahari23}\\
Triassociative, $\C$ & 1, 2 & complete list & \cite{Mainellis22tri}\\
\bottomrule
\end{tabular}
\end{table}

\subsection{Dimensions two and three}
Basri and Rikhsiboev \cite{BasriRikhsiboev07} obtained the first list of two-dimensional complex diassociative algebras by the direct method. Rakhimov later classified two-dimensional associative and diassociative algebras over an arbitrary field \cite{Rakhimov23}. In characteristic different from $2$ and $3$ his list contains eight entries, three of them with a parameter. The four algebras of \cite{BasriRikhsiboev07} correspond to four of these entries, and the other four do not appear in \cite{BasriRikhsiboev07}, so the list of \cite{BasriRikhsiboev07} has to be completed accordingly \cite{Rakhimov23}. Three-dimensional complex diassociative algebras were classified in \cite{RRB10} from the list of three-dimensional associative algebras \cite{RRB09}. In view of the correction in dimension two, an independent check of the three-dimensional list, for example by computer algebra or by the central-extension method, would be valuable.

\subsection{Nilpotent algebras}
The nilpotent diassociative algebras of dimension at most three were classified in \cite{Rikhsiboev07}, and those of dimension four in \cite{BRR15}, both by the direct method. Rakhimov \cite{Rakhimov16a} set up central extensions of associative dialgebras and applied the Skjelbred--Sund method to low-dimensional nilpotent algebras; this method gives an independent route to these lists. In dimension three the nilpotent diassociative algebras also appear inside the classification of nilpotent left-symmetric dialgebras \cite{RikhsiboevLSDA}, since every diassociative algebra is an LSDA. The two lists can be compared directly. Mainellis \cite{Mainellis22,Mainellis23} studied multipliers, covers and unicentral nilpotent diassociative algebras, extending results for Lie and Leibniz algebras; these invariants are useful for organising central extensions.

\subsection{Structure theory}
The quotient of a diassociative algebra by the ideal generated by the elements $x\lt y-x\rt y$ is an associative algebra, which links the structure of dialgebras to classical associative theory. Gonz\'alez \cite{Gonzalez13} studied associative dialgebras from this structural viewpoint. Split dialgebras, annihilator ideals and regular elements were studied by Vel\'asquez and Felipe \cite{VelasquezFelipe}, and dialgebras with a bar-unit, which are the closest analogues of unital algebras, appear in \cite{Felipe03,FLOV05,Pozhidaev08b}. Crossed modules of diassociative algebras were studied in \cite{CCKL17}, Gr\"obner--Shirshov bases in \cite{BokutChenLiu}, and the Gelfand--Kirillov dimension in \cite{ZhangChenYu}. A general construction of dialgebras from bimodules over algebras, equivalent to the Kolesnikov--Pozhidaev construction, is given in \cite{SVW16}.

\subsection{Derivations}
Derivations, which control the solvable algebras with a given nilradical and are used as invariants, have been computed for low-dimensional diassociative algebras in \cite{RRB14}; generalised $(\alpha,\beta,\gamma)$-derivations in \cite{FRSB16}; outer derivations in \cite{MQM25}; and derivations and bi-derivations in \cite{RestrepoSanchez}. For the polynomial dialgebra $\K[x,y]$ Lin and Zhang determined the derivations and automorphisms \cite{LinZhang}.

\subsection{Variants}
Hom-associative dialgebras of dimension at most four were classified by Zahari \cite{Zahari23}, and BiHom-associative dialgebras were introduced in \cite{ZahariBakayoko}. Loday and Ronco's triassociative algebras \cite{LodayRonco04} were classified in dimensions one and two by Mainellis \cite{Mainellis22tri}. Splittings of the operations of diassociative and triassociative algebras are studied in \cite{Teng25}, and Poisson dialgebras in \cite{DMM23}.

\section{Dendriform algebras}\label{sec:dend}

\begin{table}[ht]
\caption{Classification results for dendriform and related algebras.}\label{tab:dend}
\small
\begin{tabular}{@{}>{\raggedright}p{3.1cm}>{\raggedright}p{1.7cm}>{\raggedright}p{5.2cm}l@{}}
\toprule
Class & Dimension & Result & Reference\\
\midrule
Dendriform, $\C$ & 2 & $12$ classes, one with a parameter; incomplete & \cite{RRB10dend}\\
Dendriform, any field & 2 & classification for all characteristics; revises the complex list & \cite{RakhimovBekbaev25}\\
Dendriform, $\C$ & 2 & final list: one family and $13$ algebras (Theorem~\ref{thm:dend2}) & \cite{RRB10dend,RakhimovBekbaev25}\\
Zinbiel, $\C$ & 2 & complete list & \cite{Omirov02}\\
Zinbiel, $\C$ & 3 & complete list; nilpotency & \cite{DzhumadildaevTulenbaev}\\
Zinbiel, $\C$ & 4 & complete list (central extensions) & \cite{AMK23}\\
Zinbiel, $\C$ & 4 & degenerations, irreducible components & \cite{KPPV18}\\
Zinbiel, $\C$ & 5 & algebraic and geometric classification & \cite{AFK22}\\
Zinbiel, $\C$ & any & naturally graded, filiform, quasi-filiform & \cite{AdashevOmirovKhudoyberdiyev10,AdashevKhudoyberdiyevOmirov,AdashevCamachoGomezVidalKarimdjanov}\\
Tridendriform, $\C$ & low & classification & \cite{MAZ23}\\
Anti-dendriform, $\C$ & 3 & complete list & \cite{AANS}\\
\bottomrule
\end{tabular}
\end{table}

\subsection{Dimension two}
Two-dimensional complex dendriform algebras were described in \cite{RRB10dend} by the direct method: the associative algebra $(E,\ast)$ with $x\ast y=x\prec y+x\succ y$ is fixed first, and the structure constants of $\prec$ are then found. Algebras in which one of the operations vanishes are excluded there, since the other operation is then simply associative. The list of \cite{RRB10dend} contains twelve classes, one of which depends on a parameter, and it has been used to classify derivations and Rota-type operators (Rota--Baxter, Reynolds, Nijenhuis and averaging operators) on two-dimensional dendriform algebras \cite{BasdouriMosbahi24}. Rakhimov and Bekbaev \cite{RakhimovBekbaev25} classified two-dimensional dendriform algebras over an arbitrary base field, fixing instead the operation $\prec$ from the classification of all two-dimensional algebras. Over fields of characteristic different from $2$ and $3$ they obtain fourteen labelled classes, and over fields of characteristic $2$ and $3$ fifteen; they also revise the complex list of \cite{RRB10dend}. Comparing the two lists, and checking the dendriform identities and the isomorphisms by computer, gives the following final list over $\C$.

\begin{theorem}[\cite{RRB10dend,RakhimovBekbaev25}]\label{thm:dend2}
Every two-dimensional complex dendriform algebra in which both operations are non-zero is isomorphic to exactly one of the following algebras (only non-zero products are listed):
\begin{align*}
&\mathrm{Dend}_2^1(\alpha),\ \alpha\neq0: && e_1\prec e_1=e_2,\ e_1\succ e_1=\alpha e_2;\\
&\mathrm{Dend}_2^2: && e_1\prec e_1=e_1,\ e_1\succ e_2=e_2;\\
&\mathrm{Dend}_2^3: && e_1\prec e_1=e_1,\ e_2\prec e_1=e_2,\ e_1\succ e_2=e_2,\ e_2\succ e_1=-e_2;\\
&\mathrm{Dend}_2^4: && e_2\prec e_1=e_2,\ e_1\succ e_1=e_1;\\
&\mathrm{Dend}_2^5: && e_1\prec e_2=-e_2,\ e_2\prec e_1=e_2,\ e_1\succ e_1=e_1,\ e_1\succ e_2=e_2;\\
&\mathrm{Dend}_2^6: && e_1\prec e_1=e_1,\ e_2\succ e_2=e_2;\\
&\mathrm{Dend}_2^7: && e_1\prec e_2=-e_2,\ e_2\prec e_2=e_2,\ e_1\succ e_1=e_1,\ e_1\succ e_2=e_2;\\
&\mathrm{Dend}_2^8: && e_1\prec e_1=e_1+e_2,\ e_1\prec e_2=-e_2,\ e_2\prec e_2=e_2,\\ & && e_1\succ e_1=-e_2,\ e_1\succ e_2=e_2;\\
&\mathrm{Dend}_2^9: && e_1\prec e_1=e_1,\ e_2\prec e_1=e_2,\ e_2\succ e_1=-e_2,\ e_2\succ e_2=e_2;\\
&\mathrm{Dend}_2^{10}: && e_1\prec e_1=-e_2,\ e_2\prec e_1=e_2,\\ & && e_1\succ e_1=e_1+e_2,\ e_2\succ e_1=-e_2,\ e_2\succ e_2=e_2;\\
&\mathrm{Dend}_2^{11}: && e_2\prec e_1=e_2,\ e_1\succ e_1=e_1,\ e_1\succ e_2=e_2;\\
&\mathrm{Dend}_2^{12}: && e_1\prec e_1=e_1,\ e_2\prec e_1=e_2,\ e_1\succ e_2=e_2;\\
&\mathrm{Dend}_2^{13}: && e_1\prec e_2=-e_2,\ e_1\succ e_1=e_1,\ e_1\succ e_2=e_2;\\
&\mathrm{Dend}_2^{14}: && e_1\prec e_1=e_1,\ e_2\prec e_1=e_2,\ e_2\succ e_1=-e_2.
\end{align*}
Here $\mathrm{Dend}_2^1(\alpha)\cong\mathrm{Dend}_2^1(\alpha')$ if and only if $\alpha=\alpha'$. In total there are fourteen classes, one of them a one-parameter family.
\end{theorem}

The algebras $\mathrm{Dend}_2^1,\dots,\mathrm{Dend}_2^{12}$ are those of \cite{RRB10dend}, and $\mathrm{Dend}_2^{13}$, $\mathrm{Dend}_2^{14}$ are the two algebras found in \cite{RakhimovBekbaev25}. Table~\ref{tab:dend2} gives the associated associative algebra and the corresponding algebras of \cite[Theorem~3.1]{RakhimovBekbaev25}, written $\mathbb{D}$ there.

\begin{table}[ht]
\caption{Two-dimensional complex dendriform algebras with both operations non-zero.}\label{tab:dend2}
\small
\begin{tabular}{@{}lll@{}}
\toprule
Theorem~\ref{thm:dend2} & $(E,\ast)$ & \cite[Theorem~3.1]{RakhimovBekbaev25}\\
\midrule
$\mathrm{Dend}_2^1(\alpha)$ & $e_1e_1=(1+\alpha)e_2$ & $\mathbb{D}_{13}^1(\delta_1)$, $\delta_1=\alpha$\\
$\mathrm{Dend}_2^2$, $\mathrm{Dend}_2^3$ & $e_1e_1=e_1,\ e_1e_2=e_2$ & $\mathbb{D}_3^8$, $\mathbb{D}_3^{10}(0)$\\
$\mathrm{Dend}_2^4$, $\mathrm{Dend}_2^5$ & $e_1e_1=e_1,\ e_2e_1=e_2$ & $\mathbb{D}_3^4$, $\mathbb{D}_3^6$\\
$\mathrm{Dend}_2^6,\dots,\mathrm{Dend}_2^{10}$ & $\C\oplus\C$ & $\mathbb{D}_3^7$, $\mathbb{D}_1^{13}$, $\mathbb{D}_1^{14}$, $\mathbb{D}_3^9(1)\cong\mathbb{D}_3^{11}(1)$, $\mathbb{D}_3^5(1)$\\
$\mathrm{Dend}_2^{11}$, $\mathrm{Dend}_2^{12}$ & $e_1e_1=e_1,\ e_1e_2=e_2e_1=e_2$ & $\mathbb{D}_3^5(0)$, $\mathbb{D}_3^9(0)$\\
$\mathrm{Dend}_2^{13}$, $\mathrm{Dend}_2^{14}$ & $e_1e_1=e_1$ & $\mathbb{D}_7^2\cong\mathbb{D}_6^3$, $\mathbb{D}_3^{11}(0)$\\
\bottomrule
\end{tabular}
\end{table}

\begin{remark}\label{rem:dend2}
The comparison shows the following.
\begin{enumerate}
\item $\mathrm{Dend}_2^{13}$ and $\mathrm{Dend}_2^{14}$ are missing from \cite{RRB10dend} because the list of two-dimensional associative algebras used there omits $\C\oplus0$ ($e_1e_1=e_1$), which is the associative algebra of both. Moreover, $\mathrm{Dend}_2^1(0)$ of \cite{RRB10dend} has $\succ=0$ and should be excluded under the convention of that paper.
\item $\mathrm{Dend}_2^4$ is a dendriform algebra and coincides with $\mathbb{D}_3^4$ of \cite{RakhimovBekbaev25}. In \cite[Remark~3.1]{RakhimovBekbaev25} it is reproduced with $e_1\prec e_2=e_2$ in place of $e_2\prec e_1=e_2$; the resulting structure is not dendriform, which explains the statement there that $\mathrm{Dend}_2^4$ is not a dendriform algebra.
\item In \cite[Theorem~3.1]{RakhimovBekbaev25} as printed, $\mathbb{D}_7^2\cong\mathbb{D}_6^3$ via $e_1\mapsto e_1-e_2$, $e_2\mapsto-2e_2$, and $\mathbb{D}_3^{11}(1)\cong\mathbb{D}_3^9(1)$ via $e_1\mapsto e_1-e_2$, $e_2\mapsto e_2$. The structures $\mathbb{D}_3^{10}(1)$ and $\mathbb{D}_3^{12}$ do not satisfy the dendriform identities, since $x\prec y+x\succ y$ is not associative for them; presumably these are misprints.
\item If one operation vanishes, the algebra is $A^{\prec}=(A,\cdot,0)$ or $A^{\succ}=(A,0,\cdot)$ for a two-dimensional associative algebra $(A,\cdot)$. This gives twelve further non-zero algebras, two for each of the six non-zero two-dimensional associative algebras \cite{Rakhimov23}.
\end{enumerate}
\end{remark}

\subsection{Zinbiel algebras}
For Zinbiel algebras the classification is much further advanced. Two reasons are their Koszul duality with Leibniz algebras, whose well-developed classification techniques carry over, and the simplicity of their structure: a single product satisfying a single identity. Omirov classified the two-dimensional complex Zinbiel algebras \cite{Omirov02}, and Dzhumadildaev and Tulenbaev proved that finite-dimensional Zinbiel algebras over a field of characteristic zero are nilpotent and classified those of dimension three \cite{DzhumadildaevTulenbaev}; nilpotency over arbitrary fields is due to Towers \cite{Towers}. Using central extensions, Alvarez, de Mello and Kaygorodov classified the four-dimensional Zinbiel algebras and the five-dimensional ones with two-dimensional annihilator \cite{AMK23}. Alvarez, Fehlberg J\'unior and Kaygorodov completed the algebraic and geometric classification in dimension five: the variety of five-dimensional Zinbiel algebras has dimension $24$, with $16$ irreducible components and $11$ rigid algebras \cite{AFK22}. The degenerations of four-dimensional Zinbiel algebras are described in \cite{KPPV18}. In arbitrary dimension, naturally graded, filiform and quasi-filiform Zinbiel algebras are classified through the characteristic sequence \cite{AdashevOmirovKhudoyberdiyev10,AdashevKhudoyberdiyevOmirov,AdashevCamachoGomezVidalKarimdjanov}.

\subsection{Related algebras}
Tridendriform algebras of small dimension are classified in \cite{MAZ23}, and three-dimensional anti-dendriform algebras in \cite{AANS}. The connection with Rota--Baxter operators \cite{Aguiar00,EbrahimiFard02,EbrahimiFardGuo08} and $\mathcal{O}$-operators \cite{BaiGuoNi} gives a second way to produce dendriform structures: every dendriform algebra arises from an $\mathcal{O}$-operator on its associated associative algebra \cite{BaiGuoNi}, so classifications of $\mathcal{O}$-operators on small associative algebras yield dendriform algebras. New identities in dendriform algebras are studied in \cite{EMP08}, quadri-algebras in \cite{AguiarLoday04}, and dendriform algebras attached to an arbitrary variety in \cite{GubarevKolesnikov}.

\section{Cohomology and deformations}\label{sec:coh}

The classification methods of Section~\ref{sec:methods} rely on cohomology: central extensions are classified by $H^2$, and degenerations are the geometric counterpart of deformations. We summarise the relevant results in chronological order.

Loday introduced the bar homology $HY$ of diassociative algebras \cite{LodayEtAl}. Frabetti showed that for a unital associative algebra the dialgebra homology is isomorphic to Hochschild homology up to a shift of degree \cite{Frabetti97}, computed the Leibniz homology of matrix dialgebras in terms of a shuffle-type dialgebra homology \cite{Frabetti98}, and developed dialgebra (co)homology with coefficients, proving that $HY^2$ classifies abelian extensions and that Morita invariance holds for bar-unital dialgebras \cite{Frabetti01,Frabetti01b}. Majumdar and Mukherjee showed that the dialgebra cohomology $HY^*(D,D)$ carries a $G$-algebra structure, as Hochschild cohomology does \cite{MajumdarMukherjee03}, and developed the deformation theory of dialgebras controlled by this cohomology \cite{MajumdarMukherjee05}. Yau studied deformations of morphisms of dialgebras and their rigidity \cite{Yau08}, and Fialowski and Majumdar constructed miniversal deformations \cite{FialowskiMajumdar}. Later work treats equivariant \cite{Saha19} and oriented \cite{KoamSaha} dialgebras, and the varieties of dialgebras in the sense of Kolesnikov \cite{Kolesnikov08}; Bremner's survey \cite{Bremner12} covers their polynomial identities.

For the classification problem the most useful facts are these: $H^2(A,\C)$ is computable by linear algebra from the structure constants of $A$; its dimension bounds the number of parameters of the central extensions; and a vanishing second cohomology with coefficients in the algebra itself implies rigidity, so the rigid algebras that appear in a geometric classification can often be predicted from cohomology.

\section{Open problems}\label{sec:open}

We end with problems that, as far as we know, are open.

\begin{problem}
Verify the classification of three-dimensional complex diassociative algebras \cite{RRB10} independently, by computer algebra or by central extensions, in view of the correction of the two-dimensional list \cite{Rakhimov23}, and extend it to arbitrary fields as in \cite{Rakhimov23}.
\end{problem}

\begin{problem}
Classify the four-dimensional (not necessarily nilpotent) complex diassociative algebras, and the five-dimensional nilpotent ones by the central-extension method \cite{Rakhimov16a,Mainellis22}.
\end{problem}

\begin{problem}
Check the classifications of two-dimensional dendriform algebras over fields of characteristic $2$ and $3$ \cite{RakhimovBekbaev25} in the light of Remark~\ref{rem:dend2}, for example by computer algebra.
\end{problem}

\begin{problem}
Give a complete and independently verified classification of the three-dimensional complex dendriform algebras, and classify the nilpotent dendriform algebras of dimension four, for instance by central extensions.
\end{problem}

\begin{problem}
Give the geometric classification of the varieties $\cat{Dias}_n$ and $\cat{Dend}_n$ for small $n$: describe the degenerations, the irreducible components and the rigid algebras, as done for Zinbiel algebras in \cite{KPPV18,AFK22}.
\end{problem}

\begin{problem}
Develop a structure theory of diassociative algebras (radical, semisimple part, Wedderburn-type decomposition) sufficient to reduce the classification to the nilpotent case, as for associative and Lie algebras; see \cite{Gonzalez13} for first steps.
\end{problem}

\begin{problem}
Study the classification of algebras in the extended Loday diagram simultaneously, using the functors of Figures~\ref{fig:loday} and~\ref{fig:ext}: for example, compare the Leibniz algebras obtained from the classified diassociative and left-symmetric dialgebras \cite{RikhsiboevLSDA} with the known lists of Leibniz algebras.
\end{problem}

\section*{Acknowledgments}
The author used an AI assistant (Claude, Anthropic) for literature searches, drafting and language editing. All statements and references were checked by the author, who takes full responsibility for the content of this paper. This research received no external funding.

\end{document}